\documentclass[11pt,a4paper,reqno]{amsart}

\usepackage[truedimen,margin=25truemm]{geometry} 
\usepackage{amsmath,amssymb,amsthm,amscd}
\usepackage{mathtools}

\title[Cardinalities of path sets]
{Proof of a cardinality conjecture for the correction terms 
in the quantum $K$-theoretic divisor axiom}
 \author[M.~Takahashi]{Mizuki Takahashi}
 \address[Mizuki Takahashi]{Doctoral Program in Mathematics, 
 Degree Programs in Pure and Applied Sciences, 
 Graduate School of Science and Technology, 
 University of Tsukuba, 
 1-1-1 Tennodai, Tsukuba, Ibaraki 305-8571, Japan.}
 \email{s2430045@u.tsukuba.ac.jp}

\keywords{quantum $K$-theory, divisor axiom, Gromov-Witten invariants,
quantum Bruhat graph, Demazure product \newline
Mathematics Subject Classification 2020:
Primary 14N35;
Secondary 05E14, 14M15, 14N15, 20F55}

\usepackage{xcolor}

\allowdisplaybreaks
\numberwithin{equation}{section}

\newcommand{\Fg}{\mathfrak{g}}

\newcommand{\Ft}{\mathfrak{t}}

\newcommand{\BZ}{\mathbb{Z}}

\newcommand{\BC}{\mathbb{C}}

\newcommand{\CO}{\mathcal{O}}
\newcommand{\CE}{\mathcal{E}}
\newcommand{\CP}{\mathcal{P}}
\newcommand{\CW}{\mathcal{W}}
\newcommand{\CS}{\mathcal{S}}

\newcommand{\vpi}{\varpi}

\newcommand{\be}{\mathbf{e}}
\newcommand{\bp}{\mathbf{p}}
\newcommand{\bq}{\mathbf{q}}

\DeclareMathOperator{\wt}{wt}
\DeclareMathOperator{\qwt}{qwt}
\DeclareMathOperator{\ed}{end}

\newcommand{\af}{\mathrm{af}}
\newcommand{\pt}{\mathrm{pt}}
\newcommand{\Hom}{\mathrm{Hom}}
\newcommand{\QBG}{\mathrm{QBG}}

\newcommand{\Lie}{\mathrm{Lie}}
\newcommand{\LP}{\mathrm{LP}}

\newcommand{\edge}[1]{ \xrightarrow{\hspace{2pt}#1\hspace{2pt}} }

\newcommand{\rr}{\Delta_{\af}}
\newcommand{\prr}{\Delta_{\af}^{+}}
\newcommand{\Waf}{W_{\af}}
\newcommand{\Iaf}{I_{\af}}

\newcommand{\mcr}[1]{\lfloor #1 \rfloor}

\newcommand{\ti}[1]{\widetilde{#1}}
\newcommand{\ha}[1]{\widehat{#1}}

\newcommand{\ud}[1]{\underline{#1}}

\newcommand{\pair}[2]{\langle #1,\,#2 \rangle}

\newcommand{\tb}[1]{\le_{#1}}
\newcommand{\dtb}[1]{\le_{#1}^{\ast}}
\newcommand{\tbmin}[3]{\min(#1W_{#2},\le_{#3}\nobreak)}
\newcommand{\tbmax}[3]{\max(#1W_{#2},\le_{#3}^{\ast}\nobreak)}

\newcommand{\J}{J}

\newcommand{\WJs}{W_{\J}}

\newcommand{\DJs}{\Delta_{\J}^{+}}
\newcommand{\DJp}{\Delta^{+} \setminus \Delta_{\J}^{+}}

\newcommand{\QJ}{Q_{\J}}
\newcommand{\QJv}{Q_{\J}^{\vee}}
\newcommand{\QJvp}{Q_{\J}^{\vee,+}}

\newcommand{\bR}[1]{\mathbf{R}_{#1}^{\lhd}}
\newcommand{\bE}[1]{\mathbf{E}_{#1}}

\newcommand{\dast}{\mathbin{\ud{\ast}}}

\theoremstyle{plain}
\newtheorem{lem}{Lemma}[section]
\newtheorem{prop}[lem]{Proposition}
\newtheorem{thm}[lem]{Theorem}
\newtheorem{cor}[lem]{Corollary}
\newtheorem{conj}[lem]{Conjecture}

\newtheorem{ithm}{Theorem}

\theoremstyle{definition}
\newtheorem{dfn}[lem]{Definition}

\newtheorem{idfn}[ithm]{Definition}

\theoremstyle{remark}

\newtheorem{rem}[lem]{Remark}

\newenvironment{enu}{%
 \begin{enumerate}%
}{\end{enumerate}}

\begin{document}

%
\begin{abstract}
We prove a conjecture due to Kato and Sagaki \cite{KS} concerning the cardinality of sets of
directed paths in the quantum Bruhat graph that occur in the quantum
$K$-theoretic divisor axiom for the full flag manifolds in the quasi-minuscule case.
We prove that the cadinality is equal to one or even, and
its signed sum vanishes in the latter case.
Consequently, the cancellation-free divisor axiom obtained by Kato and Sagaki
under this conjecture holds unconditionally.
The proof uses Schremmer's results on affine Bruhat order and
Demazure products, together with a theorem that we prove
for arbitrary Coxeter groups.
\end{abstract}

\maketitle
%
%
\section{Introduction.}
\label{sec:intro}

Let $G$ be a connected, simply-connected, simple algebraic group over $\BC$,
let $T$ be a maximal torus of $G$, and let $B$ be a Borel subgroup
containing $T$.
Set $X:=G/B$.
Let $\Delta$ be the root system of $G$ relative to $T$,
with positive roots $\Delta^{+}$ and simple roots
$\{\alpha_j\}_{j\in I}$.
For $\alpha\in\Delta$, write $\alpha^{\vee}$ for its coroot
and $s_{\alpha}$ for the corresponding reflection.
We denote by $W=\langle s_j\mid j\in I\rangle$ the Weyl group,
where $s_j:=s_{\alpha_j}$, and by $\ell$ its length function.
Let $\vpi_j$, $j\in I$, be the fundamental weights.
Set
\[
\Lambda:=\bigoplus_{j\in I}\BZ\vpi_j,\qquad
Q^{\vee,+}:=\sum_{j\in I}\BZ_{\ge0}\alpha_j^{\vee}.
\]

For $u\in W$, let $\CO_u$ and $\CO^u$ denote the Schubert class
and the opposite Schubert class in the $T$-equivariant $K$-theory ring
$K_T(X)$, respectively.
For $\gamma_1,\ldots,\gamma_m\in K_T(X)$ and $d\in Q^{\vee,+}$,
write $\langle\gamma_1,\ldots,\gamma_m\rangle_d$ for the corresponding
$T$-equivariant $K$-theoretic Gromov-Witten invariant.
These invariants take values in
$K_T(\pt)\cong\BZ[\Lambda]$,
where $\BZ[\Lambda]$ is the group algebra of $\Lambda$.

Fix $i\in I$ such that $\vpi_i$ is quasi-minuscule in the sense that
$
\pair{\vpi_i}{\beta^{\vee}}\in\{0,1,2\}$
for every $\beta\in\Delta^{+}$.
Note that this condition holds for every fundamental weight in classical type.
By the quantum $K$-theoretic divisor axiom for $X$
due to Lenart, Naito, Sagaki, and Xu \cite[Theorem~3.2]{LNSX},
we have for $w,x\in W$ and $d\in Q^{\vee,+}$,
\begin{equation} \label{eq:ithp-qm}
\langle\CO^{s_i},\CO^w,\CO_x\rangle_d
=\langle\CO^w,\CO_x\rangle_d
-\be^{-\vpi_i+\frac12w\vpi_i+\frac12x\vpi_i}
\underbrace{\sum_{\bp\in\bR{w,x,d}}(-1)^{\ell(\bp)}}_{\text{the correction term}}.
\end{equation}
Here $\bR{w,x,d}$ is a finite set of directed paths in the quantum
Bruhat graph $\QBG(W)$; its precise definition is given in \eqref{eq:iR} below.
If $\bR{w,x,d}\ne\emptyset$, then
$\frac12(w\vpi_i+x\vpi_i)\in\Lambda$;
otherwise, the second term in \eqref{eq:ithp-qm}
is understood to be zero.

The correction term in \eqref{eq:ithp-qm} may contain cancellations.
Kato and Sagaki \cite{KS} studied this problem and obtained
a cancellation-free formula under a conjecture on the cardinalities
of the indexing sets $\bR{w,x,d}$.
Our purpose is to prove this conjecture.
Let us explain our result more precisely.

\begin{idfn} \label{idfn:QBG}
The quantum Bruhat graph $\QBG(W)$ has vertex set $W$.
For $x,y\in W$ and $\alpha\in\Delta^{+}$, it has an edge
$x\edge{\alpha}y$ labeled by $\alpha$ if $y=xs_{\alpha}$
and one of the following conditions holds:
\textup{(B)} $\ell(y)=\ell(x)+1$;
\textup{(Q)}
$\ell(y)=\ell(x)+1-2\pair{\rho}{\alpha^{\vee}}$,
where $\rho:=\frac12\sum_{\alpha\in\Delta^{+}}\alpha$.
Edges of type \textup{(B)} and \textup{(Q)} are called
Bruhat edges and quantum edges, respectively.
\end{idfn}

For a shortest directed path
$\bp:u=u_0\edge{\beta_1}u_1\edge{\beta_2}\cdots
\edge{\beta_r}u_r=v$ in $\QBG(W)$, set
$\ed(\bp):=v$ and $\ell(\bp):=r$.
Its quantum weight $\qwt(\bp)$ is defined to be the sum of $\beta_t^{\vee}$
over the quantum edges $u_{t-1}\edge{\beta_t}u_t$ in $\bp$.

Set $\J:=I\setminus\{i\}$,
$\WJs:=\langle s_j\mid j\in\J\rangle$, and
$\DJs:=\Delta^{+}\cap\bigoplus_{j\in\J}\BZ\alpha_j$.
Let $\lhd$ be an arbitrary reflection order on $\Delta^{+}$
satisfying the condition that
$\gamma\lhd\beta$
for all $\gamma\in\DJs$ and $\beta\in\DJp$.
A path is called label-increasing if its labels strictly increase
with respect to $\lhd$.

For $w,x\in W$ and $d\in Q^{\vee,+}$, consider the set
$\bR{w,x,d}$ introduced in \cite[(3.3)]{LNSX}.
In the quasi-minuscule case, we use its description
in \cite[(1.3)]{KS} as a set of directed paths
$\bp:w=u_0\edge{\beta_1}\cdots\edge{\beta_r}u_r=\ed(\bp)$
in $\QBG(W)$ satisfying
\begin{equation} \label{eq:iR}
\begin{cases}
r\ge0,\quad \beta_1\lhd\cdots\lhd\beta_r,\quad
\pair{\vpi_i}{\beta_t^{\vee}}=2\quad(1\le t\le r),\\[1mm]
\ed(\bp)\in x\WJs,\quad
\qwt(\bp)+\qwt(\ed(\bp)\Rightarrow x)\le d,\quad
\pair{\vpi_i}{d-\qwt(\bp)}=0.
\end{cases}
\end{equation}
Kato and Sagaki conjectured that
$\#\bR{w,x,d}\in\{1\}\sqcup2\BZ_{\ge0}$
\cite[Conjecture~3.5]{KS}.
Assuming this conjecture, they proved that
$\langle\CO^{s_i},\CO^w,\CO_x\rangle_d=1$
whenever $\#\bR{w,x,d}\ge2$ \cite[Theorem~3.7]{KS}.
They also obtained this equality without the conjecture when
the coefficients of $d$ along the simple coroots indexed by $\J$
are sufficiently large \cite[Corollary~3.11]{KS}.
Our main result removes this restriction on $d$.

\begin{ithm}[see Theorem~\ref{thm:main}] \label{ithm:main}
Let $i\in I$ be such that $\vpi_i$ is quasi-minuscule.
For every $w,x\in W$ and $d\in Q^{\vee,+}$, we have
\begin{equation} \label{eq:main11}
\#\bR{w,x,d}\in\{1\}\sqcup2\BZ_{\ge0}.
\end{equation} 
Moreover, if $\#\bR{w,x,d}\ge2$, then
\begin{equation} \label{eq:main21}
\sum_{\bp\in\bR{w,x,d}}(-1)^{\ell(\bp)}=0.
\end{equation}
\end{ithm}

The cardinality statement includes the possibility that
$\bR{w,x,d}$ is empty.
We prove (\ref{eq:main11}) and (\ref{eq:main21}) simultaneously,
without using \cite[Theorem~3.7]{KS}.

This paper is organized as follows.
In Section~\ref{sec:cox}, we prove Theorem~\ref{thm:cox}
on cardinalities and signed sums of sets defined by the Demazure product
in an arbitrary Coxeter group.
In section~\ref{sec:prelim},we recalls the notation for root systems,
quantum Bruhat graphs, and tilted Bruhat orders, and establish
the lemmas needed later.
In section~\ref{sec:main}, we state the main results (Theorem~\ref{thm:main}),
and then reduce its proof
to the study of the set of "end points''.
In Section~\ref{sec:aff}, we apply the results mentioned in Section~\ref{sec:cox}
to the affine Weyl group
and derive the formulas used in Section~\ref{sec:proof} to construct a bijection.
Finally, in Section~\ref{sec:prfmain}, we complete the proof of the main theorem.

%
\subsection*{Acknowledgements}

The author would like to thank Professor Daisuke Sagaki
for his guidance and helpful advice.
The author also thanks Professor Tomoyuki Kakehi
for his guidance during the author's master's studies.
This work was supported by JST SPRING,
Japan Grant Number JPMJSP2124.
%
%
\section{Demazure products and Bruhat intervals in Coxeter groups.}
\label{sec:cox}

Throughout this section, let $(\CW,\CS)$ be an arbitrary Coxeter system.
We denote by $\ell : \CW \rightarrow \BZ_{\ge 0}$ the length function,
by $e \in \CW$ the identity element, and by $\le$ the Bruhat order on $\CW$.

%
\subsection{Bruhat order and Demazure products.}
\label{subsec:Demazure}
\begin{lem}[{Lifting property;
see, e.g., \cite[Proposition~2.2.7]{BB}}]
\label{lem:lifting}
Let $f,g\in\CW$ and $s\in\CS$ satisfy
$f<g$, $fs>f$, and $gs<g$.
Then
\begin{equation} \label{eq:lifting}
f\le gs \qquad\text{and}\qquad fs\le g.
\end{equation}
\end{lem}

We denote by $\ast$ the Demazure product on $\CW$,
that is, the unique associative product defined by
$f \ast s = \max\{f,fs\}$ (with respect to the Bruhat order $\leq$)
for $f \in \CW$ and $s \in \CS$; see \cite[Section~3]{BM}.
For $g\in\CW$, choose a reduced expression
$g=s_{j_1}s_{j_2}\cdots s_{j_n}$, and set
\[
f\ast g
:=\bigl(\cdots((f\ast s_{j_1})\ast s_{j_2})
\cdots\bigr)\ast s_{j_n}.
\]
This definition is independent of the choice of a reduced
expression.
%
%
\begin{prop}[{see \cite[Propositions~3.1 and 3.2]{BM}}]
\label{prop:He}
Let $f,g \in \CW$. Then, $f \ast g$ is equal to the unique maximal element of the subset
$\{ f'g' \mid f' \le f, \, g' \le g \}$.
Moreover, there exists $f' \in \CW$ such that
$f' \le f$, $f \ast g = f'g$, and $\ell(f \ast g) = \ell(f') + \ell(g)$.
\end{prop}

\begin{proof}
For $a,b \in \CW$ such that $a \le f$ and $b \le g$,
we have
\begin{equation*}
ab \le a \ast b \le f \ast g
\end{equation*}
by \cite[Proposition~3.1(c), (d)]{BM}.
Set $f' := (f \ast g)g^{-1}$.
By \cite[Proposition~3.1(e)]{BM}, we have
$f' \le f$ and $f'g = f' \ast g = f \ast g$.
Hence $f \ast g$ belongs to the given subset
and is its unique maximal element.
Finally, since $f' \ast g = f'g$, we obtain
$\ell(f \ast g) = \ell(f') + \ell(g)$
from \cite[Proposition~3.2]{BM}.
\end{proof}

%
%

For $f \in \CW$ and $s \in \CS$, we define
$f \dast s := \min\{f,fs\}$ (with respect to the Bruhat order).
For $g \in \CW$, choose a reduced expression
$g=s_{j_1}s_{j_2}\cdots s_{j_n}$, and set
\[
f \dast g
:= \bigl(\cdots((f \dast s_{j_1}) \dast s_{j_2})
\cdots\bigr) \dast s_{j_n}.
\]

\begin{lem} \label{lem:dast}
Let $f,g \in \CW$.
The definition of $f \dast g$ above is independent
of the choice of a reduced expression for $g$.
Moreover, $f \dast g$ is equal to the unique minimal element of
$M(f,g) := \{fv \mid v \in \CW,\ v \le g\}$
with respect to the Bruhat order.
In addition, the following hold.
\begin{enu}
\item If $s \in \CS$ satisfies $gs > g$, then
$f \dast (gs)=\min\{f \dast g,\,(f \dast g)s\}$.

\item If $s \in \CS$ satisfies $sg > g$, then
$f \dast (sg)=\min\{f,fs\}\dast g$.

\item If $g' \in \CW$ satisfies $g' \le g$, then
$f \dast g \le f \dast g'$.
\end{enu}
\end{lem}
\begin{proof}
First, for $a,b \in \CW$ with $a \le b$ and $s \in \CS$, we have
$
\min\{a,as\} \le \min\{b,bs\}.
$
Indeed, if $bs>b$, then $\min\{a,as\}\le a\le b=\min\{b,bs\}$.
If $bs<b$, set $m:=\min\{a,as\}$.
Since $ms>m$ and $m<b$, the lifting property gives
$m\le bs=\min\{b,bs\}$.

Let $g = s_{j_1}s_{j_2} \cdots s_{j_n}$ be a reduced expression for $g$,
and define $m_{0},m_{1},\dots,m_{n} \in \CW$ by
$m_{0} := f$, and
$m_{t} := \min \{ m_{t-1}, \, m_{t-1}s_{j_t} \}$
for $1 \le t \le n$.
We claim that
\begin{equation} \label{eq:dast1}
\text{$m_{n} \in M(f,g)$, and
$m_{n} \le h$ for all $h \in M(f,g)$}.
\end{equation}
We show \eqref{eq:dast1} by induction on $n = \ell(g)$.
If $n = 0$, then the assertion is obvious since
$M(f,e) = \{f\}$.
Assume that $n > 0$, and set
$g' := s_{j_1} \cdots s_{j_{n-1}}$, $s := s_{j_n}$.
By the induction hypothesis, we have
$m_{n-1} \in M(f,g')$,
and $m_{n-1} \le h$ for all $h \in M(f,g')$.
Also, since $g's = g > g'$, we see by
the subword property (see \cite[Theorem~2.2.2]{BB}) that
$
M(f,g)
=
M(f,g') \cup \{ hs \mid h \in M(f,g') \}.
$
Hence 
$m_{n} \in \{ m_{n-1}, \, m_{n-1}s \} \subset M(f,g)$.
Let $h \in M(f,g')$.
Since $m_{n-1} \le h$, we see by
the observation above that
$
m_{n}
=
\min \{ m_{n-1}, \, m_{n-1}s \}
\le \min \{h,hs\},
$
which implies that $m_{n} \le h$ and $m_{n} \le hs$.
Thus we have shown \eqref{eq:dast1}.

By \eqref{eq:dast1}, the set $M(f,g)$ has
a unique minimal element $f \dast g = m_{n}$;
in particular, $m_{n}$ does not depend on the choice
of a reduced expression for $g$.
Part~(1) follows by the construction to
above applied a reduced expression for $gs$ ending in $s$.
Similarly, part~(2) follows by the construction above applied to
a reduced expression for $sg$ beginning with $s$.
Part~(3) follows from \eqref{eq:dast1}, since
$M(f,g') \subset M(f,g)$ if $g' \le g$.
Thus we have proved the lemma.
\end{proof}
%
%
\subsection{Cardinalities and signed sums.}
\label{subsec:cardinalities}

We fix $u \in \CW$. For $f,g \in \CW$, we set
%
%
\begin{equation} \label{eq:CE}
\CE_{u}(f,g) := \bigl\{ h \in \CW \mid h \le g, \ u \ast h = f \bigr\},
\end{equation}
which is a finite set.

%
%
\begin{prop} \label{prop:decomp}
Let $f,g \in \CW$, and let $s \in \CS$ be such that $gs < g$.
\begin{enu}
\item If $fs > f$, then $\CE_{u}(f,g) = \CE_{u}(f,gs)$.

\item Assume that $fs < f$, and set
\begin{equation} \label{eq:CP}
\CP(f,g\,;s) := \bigl\{ h \in \CE_{u}(f,g) \mid hs \in \CE_{u}(f,g) \bigr\}.
\end{equation}
Then, the map $h \mapsto hs$ is a fixed-point-free involution on $\CP(f,g\,;s)$.
Moreover, the map
\begin{equation} \label{eq:decomp}
\CE_{u}(fs,gs) \rightarrow \CE_{u}(f,g) \setminus \CP(f,g\,;s), \qquad h' \mapsto h's,
\end{equation}
is bijective, and $\ell(h's) = \ell(h') + 1$ for all $h' \in \CE_{u}(fs,gs)$.
\end{enu}
\end{prop}

\begin{proof}
(1) Let $h \in \CE_{u}(f,g)$.
Suppose, for a contradiction, that $hs < h$.
Since $h \ast s = h$, it follows by
the associativity of $\ast$ that
\begin{equation*}
f \ast s
= (u \ast h) \ast s
= u \ast (h \ast s)
= u \ast h
= f.
\end{equation*}
This contradicts $f \ast s = fs > f$.
Thus we get $hs > h$.
Since $h \le g$ and $gs < g$, it follows from the lifting property that
$h \le gs$, and hence $h \in \CE_{u}(f,gs)$.
The opposite inclusion $\CE_{u}(f,g) \supset \CE_{u}(f,gs)$ is obvious since $gs < g$.

(2) It is obvious from the definition that
the map $h \mapsto hs$ is a fixed-point-free involution on $\CP(f,g\,;s)$.
We show that $h's \in \CE_{u}(f,g)\setminus \CP(f,g\,;s)$ for all $h' \in \CE_{u}(fs,gs)$.
Let $h' \in \CE_{u}(fs,gs)$.
Since $(fs)s = f > fs$,
the same argument as in part (1) shows that $h's > h'$.
Since $h' \le gs$, we have $h's \le g$ by the lifting property.
Also, we obtain
\begin{equation*}
u \ast (h's)
= u \ast (h' \ast s)
= (u \ast h') \ast s
= (fs) \ast s
= f.
\end{equation*}
Hence we get $h's \in \CE_{u}(f,g)$.
Since $u \ast h' = fs \ne f$, we have
$(h's)s = h' \notin \CE_{u}(f,g)$, and hence $h's \notin \CP(f,g\,;s)$.
Also, we see that
$\ell(h's) = \ell(h') + 1$ for all $h' \in \CE_{u}(fs,gs)$,
and that the map \eqref{eq:decomp} is injective.

It remains to show that the map \eqref{eq:decomp} is surjective.
Let $h \in \CE_{u}(f,g) \setminus \CP(f,g\,;s)$.
Suppose, for a contradiction, that $hs > h$.
Since $h \le g$ and $gs < g$, we have $hs \le g$ by the lifting property.
Also, we obtain
\begin{equation*}
u \ast (hs)
= u \ast (h \ast s)
= (u \ast h) \ast s
= f \ast s
= f.
\end{equation*}
Hence we get $hs \in \CE_{u}(f,g)$,
which contradicts $h \notin \CP(f,g\,;s)$.
Thus we obtain $hs < h$.
Since $(hs) \ast s = h$,
\begin{equation*}
(u \ast (hs)) \ast s
= u \ast ((hs) \ast s)
= u \ast h
= f.
\end{equation*}
It follows that
$f \in \{u \ast (hs),\, (u \ast (hs))s\}$,
and hence $u \ast (hs) \in \{f,fs\}$.
If $u \ast (hs) = f$, then $hs \in \CE_{u}(f,g)$ since $hs < h \le g$,
which contradicts $h \notin \CP(f,g\,;s)$.
Therefore, we get $u \ast (hs) = fs$.
Since $hs \le gs$ by the lifting property,
we conclude that $hs \in \CE_{u}(fs,gs)$.
This proves the surjectivity of the map \eqref{eq:decomp}.
\end{proof}

We are now ready to state and prove the main result of this section.

\begin{thm} \label{thm:cox}
Let $u,f,g \in \CW$.
\begin{enu}
\item If $f \dast g^{-1} = u$, then
$\CE_{u}(f,g) = \{u^{-1}f\}$.

\item If $f \dast g^{-1} \ne u$, then
$\#\CE_{u}(f,g)$ is an even number
{\rm(}possibly zero{\rm)}.
\end{enu}
Moreover, we have
\begin{equation} \label{eq:cox3}
\sum_{h \in \CE_{u}(f,g)} (-1)^{\ell(h)}
=
\begin{cases}
(-1)^{\ell(f)-\ell(u)}
  & \text{\rm if $f \dast g^{-1} = u$}, \\[1mm]
0 & \text{\rm otherwise}.
\end{cases}
\end{equation}
\end{thm}

\begin{proof}
We first prove (1), (2), and the vanishing in \eqref{eq:cox3}
when $f \dast g^{-1} \ne u$,
simultaneously for all $f \in \CW$ by induction on $\ell(g)$.
Assume that $\ell(g) = 0$, i.e., $g = e$.
Since $u \ast e = u$ and $f \dast e = f$, we see that
if $f \dast e = u$, then $\CE_{u}(f,e) = \{e\} = \{ u^{-1}f \}$, and that
if $f \dast e \ne u$, then $\CE_{u}(f,e) = \emptyset$.
Hence the assertions hold in this case.

Assume that $\ell(g) > 0$, and take $s \in \CS$ such that $gs < g$.
Since $g^{-1} = s(gs)^{-1} > (gs)^{-1}$, it follows from Lemma~\ref{lem:dast}\,(2) that
\begin{equation} \label{eq:cox1}
f \dast g^{-1} = \min \{f,fs\} \dast (gs)^{-1} =
\begin{cases}
f \dast (gs)^{-1} & \text{if $fs > f$}, \\[1mm]
(fs) \dast (gs)^{-1} & \text{if $fs < f$}.
\end{cases}
\end{equation}

If $fs > f$, then we have $\CE_{u}(f,g) = \CE_{u}(f,gs)$ by Proposition~\ref{prop:decomp}\,(1), and
$f \dast g^{-1} = f \dast (gs)^{-1}$ by \eqref{eq:cox1}.
Hence the assertions follow from the induction hypothesis.

If $fs < f$, then Proposition~\ref{prop:decomp}\,(2) gives
\begin{equation} \label{eq:cox2}
\CE_{u}(f,g) = \CP(f,g\,;s) \sqcup \bigl\{ h's \mid h' \in \CE_{u}(fs,gs) \bigr\},
\end{equation}
while \eqref{eq:cox1} gives $f \dast g^{-1} = (fs) \dast (gs)^{-1}$.

Assume first that $f \dast g^{-1} = u$. Since $(fs) \dast (gs)^{-1} = u$,
it follows from the induction hypothesis that
$\CE_{u}(fs,gs) = \{ u^{-1}fs \}$.
Hence we see by \eqref{eq:cox2} that $u^{-1}f = (u^{-1}fs)s \in \CE_{u}(f,g)$.
Conversely, let $h \in \CE_{u}(f,g)$.
By Proposition~\ref{prop:He}, there exists $a \leq u$ such that $f = ah$,
and hence $fh^{-1} = a \leq u$.
Since $h^{-1} \leq g^{-1}$, the characterization in
Lemma~\ref{lem:dast} gives
$
    u = f \mathbin{\underline{\ast}} g^{-1}
    \leq fh^{-1} \leq u.
$
Hence $fh^{-1} = u$, and therefore $h = u^{-1}f$.
Together with $u^{-1}f \in \CE_{u}(f,g)$,
this gives $\CE_{u}(f,g) = \{u^{-1}f\}$, as desired.

Assume next that $f \dast g^{-1} \ne u$. Since $(fs) \dast (gs)^{-1} \ne u$,
it follows from the induction hypothesis that
$\# \CE_{u}(fs,gs)$ is an even number, and $\sum_{h' \in \CE_{u}(fs,gs)} (-1)^{\ell(h')} = 0$.
Also, since the map $h \mapsto hs$ is an involution on $\CP(f,g\,;s)$ without fixed points,
and since $(-1)^{\ell(hs)} = -(-1)^{\ell(h)}$, we see that
$\# \CP(f,g\,;s)$ is an even number, and $\sum_{h \in \CP(f,g\,;s)} (-1)^{\ell(h)} = 0$.
Therefore, we deduce from \eqref{eq:cox2},
together with the equality $\ell(h's) = \ell(h') + 1$ for $h' \in \CE_{u}(fs,gs)$, that
$\# \CE_{u}(f,g) = \# \CP(f,g\,;s) + \# \CE_{u}(fs,gs)$ is an even number, and
\begin{equation*}
\sum_{h \in \CE_{u}(f,g)} (-1)^{\ell(h)} =
\sum_{h \in \CP(f,g\,;s)} (-1)^{\ell(h)} - \sum_{h' \in \CE_{u}(fs,gs)} (-1)^{\ell(h')} = 0,
\end{equation*}
as desired. 

It remains to prove \eqref{eq:cox3} when $f \dast g^{-1} = u$.
In this case, $\CE_{u}(f,g) = \{u^{-1}f\}$ by (1).
Since $u \ast (u^{-1}f) = f$,
Proposition~\ref{prop:He} gives an element $f' \in \CW$ such that
\[
f = f'(u^{-1}f)
\qquad \text{and} \qquad
\ell(f) = \ell(f') + \ell(u^{-1}f).
\]
The first equality implies $f' = u$, and hence
$\ell(u^{-1}f) = \ell(f) - \ell(u)$.
Thus \eqref{eq:cox3} follows.
This completes the proof of Theorem~\ref{thm:cox}.
\end{proof}

%
%

\section{Root systems and quantum Bruhat graphs.}
\label{sec:prelim}
%
%
\subsection{Notation for root systems.} \label{subsec:alggrp}
In this subsection, we recall the basic notation
from \cite{KS} and \cite{LNSX}.
Let $G$ be a connected, simply-connected, simple (linear) algebraic group over $\BC$,
and $T$ a maximal torus of $G$.
We set $\Fg := \Lie(G)$ and $\Ft := \Lie(T)$;
$\Fg$ is a finite-dimensional simple Lie algebra over $\BC$, and
$\Ft$ is a Cartan subalgebra of $\Fg$.
We denote by $\pair{\cdot\,}{\cdot} : \Ft^{\ast} \times \Ft \rightarrow \BC$
the canonical pairing, where $\Ft^{\ast} := \Hom_{\BC}(\Ft, \BC)$.
Let $\Delta \subset \Ft^{\ast}$ be the root system of $\Fg$,
$\Delta^{+} \subset \Delta$ the set of positive roots,
and $\{ \alpha_{j} \}_{j \in I} \subset \Delta^{+}$ the set of simple roots.
Denote by $\Delta^{-}:=-\Delta^{+}$ the set of negative roots.
We denote by $\alpha^{\vee} \in \Ft$ the coroot of $\alpha \in \Delta$.
Also, we denote by $\theta \in \Delta^{+}$ the highest root of $\Delta$,
and set $\rho := (1/2) \sum_{\alpha \in \Delta^{+}} \alpha$.
The root lattice $Q$ and the coroot lattice $Q^{\vee}$ of $\Fg$ are defined by
$Q := \bigoplus_{j \in I} \BZ \alpha_{j}$ and
$Q^{\vee} := \bigoplus_{j \in I} \BZ \alpha_{j}^{\vee}$, respectively.
We set $Q^{+} := \sum_{j \in I} \BZ_{\ge 0} \alpha_{j}$ and
$Q^{\vee,+} := \sum_{j \in I} \BZ_{\ge 0} \alpha_{j}^{\vee}$.
For $\xi,\zeta \in \Ft$, we write $\xi \ge \zeta$ if $\xi-\zeta \in Q^{\vee,+}$.
For $i \in I$, the element $\vpi_{i} \in \Ft^{\ast}$
given by $\pair{\vpi_{i}}{\alpha_{j}^{\vee}} = \delta_{i,j}$ for all $j \in I$,
with $\delta_{i,j}$ the Kronecker delta, is called the $i$-th fundamental weight.
Denote by $\Lambda := \bigoplus_{j \in I} \BZ \vpi_{j}$
the (integral) weight lattice of $\Fg$,
and by $\Lambda^{+} := \sum_{j \in I} \BZ_{\ge 0} \vpi_{j}$
the set of dominant (integral) weights.
We denote by $\BZ[\Lambda]$ the group algebra of $\Lambda$, that is,
the associative $\BZ$-algebra with a $\BZ$-basis $\bigl\{ \be^{\nu} \mid \nu \in \Lambda \bigr\}$,
where the product is defined by $\be^{\mu} \be^{\nu} : =
\be^{\mu + \nu}$ for $\mu,\,\nu \in \Lambda$.

A reflection $s_{\alpha} \in GL(\Ft^{\ast})$, $\alpha \in \Delta$,
is defined by $s_{\alpha} \mu := \mu - \pair{\mu}{\alpha^{\vee}} \alpha$
for $\mu \in \Ft^{\ast}$. Denote by $s_{j} := s_{\alpha_{j}}$ for $j \in I$
the simple reflection in $\alpha_{j}$.
The (finite) Weyl group $W$ of $\Fg$ is defined to be the subgroup of $GL(\Ft^{\ast})$
generated by $\{ s_{j} \}_{j \in I}$, that is, $W : = \langle s_{j} \mid j \in I \rangle$.
For $w \in W$, we denote by $\ell(w)$ the length of $w$.
Denote by $e \in W$ the identity element, and by $\le$ the Bruhat order on $W$.
Recall that $W$ acts also on $\Ft$ in such a way that
$\pair{w\mu}{w\xi} = \pair{\mu}{\xi}$ for $w \in W$, $\mu \in \Ft^{\ast}$, and $\xi \in \Ft$;
we have $w\alpha^{\vee} = (w\alpha)^{\vee}$ for $w \in W$ and $\alpha \in \Delta$.

Let $\J$ be a subset of $I$. We set
\begin{align*}
& \QJ := \bigoplus_{j \in \J} \BZ \alpha_j, \qquad
  \Delta^{\pm}_{\J} := \Delta^{\pm} \cap \QJ, \qquad
  \WJs := \langle s_{j} \mid j \in \J \rangle, \\
& \QJv := \bigoplus_{j \in \J} \BZ \alpha_{j}^{\vee}, \qquad
  \QJvp := \sum_{j \in \J} \BZ_{\ge 0} \alpha_{j}^{\vee}. 
\end{align*}
For $w \in W$, let $\mcr{w}=\mcr{w}^{\J}$
denote the minimal(-length) coset representative for the coset $w\WJs$.
We know from \cite[Proposition~2.5.1]{BB} that for $w,v \in W$,
\begin{equation} \label{eq:proj}
w \le v \quad \Longrightarrow \quad \mcr{w}^{\J} \le \mcr{v}^{\J}.
\end{equation}

%
\subsection{Quantum Bruhat graph.}
\label{subsec:QBG}

\begin{dfn} \label{dfn:QBG}
The quantum Bruhat graph $\QBG(W)$ is the $\Delta^{+}$-labeled
directed graph whose vertices are the elements of $W$ and
whose edges are of the following form:
$x \edge{\alpha} y$, with $x, y \in W$ and $\alpha \in \Delta^{+}$,
such that $y = x s_{\alpha}$ and either of the following holds:
(B) $\ell(y) = \ell (x) + 1$;
(Q) $\ell(y) = \ell (x) + 1 - 2 \pair{\rho}{\alpha^{\vee}}$.
An edge satisfying (B) (resp., (Q)) is called a Bruhat edge (resp., quantum edge).
\end{dfn}

Let
\begin{equation} \label{eq:dp}
\bp: w = u_{0} \edge{\beta_{1}} u_{1} \edge{\beta_{2}} \cdots \edge{\beta_{r}} u_{r} = v
\end{equation}
be a directed path in the quantum Bruhat graph $\QBG(W)$;
we simply write $\bp$ as $\bp: w \edge{\beta_1,\dots,\beta_r} v$.
We set $\ed(\bp) := u_{r} = v$ and $\ell(\bp)=r$.
A directed path is called trivial its length is equal to $0$.
For a directed path $\bp$ of the form \eqref{eq:dp}, we set
\begin{equation*}
\qwt(\bp) :=  \sum_{ \begin{subarray}{c} 1 \le t \le r \\
\text{$u_{t-1} \edge{\beta_t} u_{t}$ is} \\
\text{a quantum edge} \end{subarray} } \beta_{t}^{\vee} \in Q^{\vee,+}.
\end{equation*}
Let $w,v \in W$. By \cite[Theorem~7.3]{LNSSS1}, there exists a directed
path from $w$ to $v$ in $\QBG(W)$. We denote by $\ell(w \Rightarrow v)$
the length of a shortest such path.
Let $\bp$ be a shortest directed path from $w$ to $v$ in $\QBG(W)$.
We set $\qwt(w \Rightarrow v) : = \qwt(\bp)$; we know from \cite[Proposition~8.1]{LNSSS1} that
$\qwt(w \Rightarrow v)$ does not depend on the choice of a shortest directed path $\bp$.
%
%
\begin{prop}[{\cite[Proposition~8.1]{LNSSS1}}] \label{prop:qwt}
Let $w,v \in W$, and let $\bq$ be
a directed path (not necessarily shortest)
from $w$ to $v$ in $\QBG(W)$. Then we have
$\qwt(\bq) \ge \qwt(w \Rightarrow v)$.
\end{prop}

Applying Proposition~\ref{prop:qwt} to the concatenation of shortest
directed paths from $w$ to $u$ and from $u$ to $v$, we obtain
\begin{equation} \label{eq:tri}
\qwt(w \Rightarrow v) \le
\qwt(w \Rightarrow u) + \qwt(u \Rightarrow v)
\end{equation}
for all $w,u,v \in W$.
For any directed path $\bp$ from $w$ to $v$ in $\QBG(W)$, we have
\[
(-1)^{\ell(\bp)} = (-1)^{\ell(v)-\ell(w)}.
\]
%
\begin{dfn}[{cf. \cite[Section~3.2]{LNSSS2}}] \label{dfn:QBa}
Let $\lambda \in \Lambda^{+}$. For a rational number $0 \le a < 1$,
we define $\QBG_{a\lambda}(W)$
to be the subgraph of $\QBG(W)$
with the same vertex set but having only those directed edges of
the form $x \edge{\alpha} y$ for which
$a\pair{\lambda}{\alpha^{\vee}} \in \BZ$ holds.
\end{dfn}

Let $\lhd$ be a reflection (convex) order on $\Delta^{+}$;
see, e.g., \cite[Section 2.2]{KNS}.
A directed path $\bp: w \edge{\beta_1,\dots,\beta_r} v$ of the form \eqref{eq:dp}
is said to be label-increasing with respect to $\lhd$ if $\beta_{1} \lhd \cdots \lhd \beta_{r}$;
the trivial directed paths are regarded as label-increasing ones.

%
%
%
\begin{lem} \label{lem:const}
Let $J$ be a subset of $I$.
Let $w \in W$, and let $v_1,v_2 \in W$ be such that $v_2 \in v_1\WJs$.
Then we have $\pair{\vpi_{k}}{\qwt(w \Rightarrow v_1)} = \pair{\vpi_{k}}{\qwt(w \Rightarrow v_2)}$
for all $k \in I \setminus \J$.
Namely, $\qwt(w \Rightarrow v_1) - \qwt(w \Rightarrow v_2) \in \QJv$.
\end{lem}

\begin{proof}
Let $k \in I \setminus \J$.
By \eqref{eq:tri}, we have $\qwt(w \Rightarrow v_2) \le \qwt(w \Rightarrow v_1) + \qwt(v_1 \Rightarrow v_2)$.
Since $\qwt(v_1 \Rightarrow v_2) \in \QJvp$ by \cite[Lemma~2.9]{KS}, 
we deduce that
$\pair{\vpi_{k}}{\qwt(w \Rightarrow v_2)} \le \pair{\vpi_{k}}{\qwt(w \Rightarrow v_1)}$.
By interchanging the roles of $v_1$ and $v_2$, we obtain the opposite inequality.
\end{proof}

%
\begin{lem} \label{lem:stay}
Let $J$ be a subset of $I$.
Let $x \in W$, and let $v' \in x\WJs$.
If $v \in W$ satisfies
$\qwt(v' \Rightarrow v) + \qwt(v \Rightarrow x) \in \QJvp$,
then $v \in x\WJs$.
\end{lem}

\begin{proof}
Let
$\bq : v' = y_{0} \edge{\gamma_{1}} y_{1} \edge{\gamma_{2}} \cdots \edge{\gamma_{m}} y_{m} = x$
be the concatenation of a shortest directed path from $v'$ to $v$ and
a shortest directed path from $v$ to $x$ in $\QBG(W)$;
note that $v = y_{t_{0}}$ for some $0 \le t_{0} \le m$, and that
$\qwt(\bq) = \qwt(v' \Rightarrow v) + \qwt(v \Rightarrow x) \in \QJvp$ by the assumption.
We claim that
\begin{equation} \label{eq:stay1}
\mcr{y_{t-1}}^{\J} \le \mcr{y_{t}}^{\J} \qquad \text{for all $1 \le t \le m$}.
\end{equation}
Indeed, if $y_{t-1} \edge{\gamma_{t}} y_{t}$ is a Bruhat edge,
then we have $y_{t-1} < y_{t}$, and hence
$\mcr{y_{t-1}}^{\J} \le \mcr{y_{t}}^{\J}$ by \eqref{eq:proj}.
Assume that $y_{t-1} \edge{\gamma_{t}} y_{t}$ is a quantum edge.
Since $\qwt(\bq)$ is the sum of $\gamma_{t}^{\vee} \in Q^{\vee,+}$
over all quantum edges in $\bq$, and since $\qwt(\bq) \in \QJvp$,
it follows that $\gamma_{t}^{\vee} \in \QJvp$, and hence
$\gamma_{t} \in \DJs$.
Therefore, we have $y_{t} = y_{t-1}s_{\gamma_{t}} \in y_{t-1}\WJs$,
which implies that $\mcr{y_{t-1}}^{\J} = \mcr{y_{t}}^{\J}$.
Thus we have shown \eqref{eq:stay1}.
Because $y_{0} = v' \in x\WJs$ and $y_{m} = x$, we see by \eqref{eq:stay1} that
\begin{equation*}
\mcr{x}^{\J} = \mcr{y_{0}}^{\J} \le \mcr{y_{1}}^{\J} \le \cdots \le \mcr{y_{m}}^{\J} = \mcr{x}^{\J},
\end{equation*}
which implies that $\mcr{y_{t}}^{\J} = \mcr{x}^{\J}$ for all $0 \le t \le m$;
in particular, $v = y_{t_{0}} \in x\WJs$.
Thus we have proved the lemma.
\end{proof}

%
\subsection{Tilted Bruhat orders.}
\label{subsec:tilted}
%
%
\begin{dfn}[{tilted Bruhat order; see \cite[Section~7.1]{LNSSS1}}] \label{dfn:tilted}
For each $w \in W$, we define the $w$-tilted Bruhat order $\tb{w}$ on $W$ as follows:
for $v_{1},v_{2} \in W$,
\begin{equation} \label{eq:tilted}
v_{1} \tb{w} v_{2} \iff \ell(w \Rightarrow v_{2}) =
 \ell(w \Rightarrow v_{1}) + \ell(v_{1} \Rightarrow v_{2}).
\end{equation}
\end{dfn}
%
%
\begin{prop}[{\cite[Theorem~7.1]{LNSSS1}}] \label{prop:tbmin}
Let $\J$ be a subset of $I$, and let $w \in W$.
Then each coset $x\WJs$ for $x \in W$ has a unique minimal element
with respect to $\tb{w}${\rm;} we denote it by $\tbmin{x}{\J}{w}$.
\end{prop}
%
%
\begin{dfn}[{dual tilted Bruhat order; see \cite{NOS}}] \label{dfn:dtilted}
For each $v \in W$, we define the dual $v$-tilted Bruhat order $\dtb{v}$ on $W$ as follows:
for $w_{1},w_{2} \in W$,
\begin{equation} \label{eq:dtilted}
w_{1} \dtb{v} w_{2} \iff \ell(w_{1} \Rightarrow v) =
 \ell(w_{1} \Rightarrow w_{2}) + \ell(w_{2} \Rightarrow v).
\end{equation}
\end{dfn}
%
%
\begin{prop}[{\cite[Proposition~2.25]{NOS}}] \label{prop:tbmax}
Let $\J$ be a subset of $I$, and let $v \in W$.
Then each coset $w\WJs$ for $w \in W$ has a unique maximal element
with respect to $\dtb{v}${\rm;} we denote it by $\tbmax{w}{\J}{v}$.
\end{prop}
%
%
Set $w' := \tbmax{w}{\J}{v}$. By Proposition~\ref{prop:tbmax}
and \eqref{eq:dtilted}, we have
\begin{equation} \label{eq:tbmax1}
\ell(u \Rightarrow v)
= \ell(u \Rightarrow w') + \ell(w' \Rightarrow v)
\qquad \text{for all $u \in w\WJs$}.
\end{equation}
Thus, for all $u \in w\WJs$,
$\qwt(u \Rightarrow v)
= \qwt(u \Rightarrow w') + \qwt(w' \Rightarrow v)$.
%
%
\section{Main results.}
\label{sec:main}

Fix $i \in I$ such that $\vpi_{i}$ is quasi-minuscule in the sense that
$\pair{\vpi_{i}}{\beta^{\vee}} \in \{0,1,2\}$ for all $\beta \in \Delta^{+}$;
for the list of quasi-minuscule fundamental weights, see \cite[(3.3)]{KS}.
Set $\J:=I \setminus \{i\}$, and
\begin{equation} \label{eq:D2}
\Delta^{+}_{2}:=\bigl\{ \beta \in \Delta^{+} \mid \pair{\vpi_{i}}{\beta^{\vee}} = 2 \bigr\};
\end{equation}
note that $\DJs = \bigl\{ \beta \in \Delta^{+} \mid \pair{\vpi_{i}}{\beta^{\vee}} = 0 \bigr\}$.
As in \cite[(2.3)]{LNSX} and \cite[(2.7)]{KS}, let (and fix) $\lhd$ be an arbitrary reflection
order on $\Delta^{+}$ satisfying the condition that
%
%
\begin{equation} \label{eq:ro}
\gamma \lhd \beta \quad
\text{for all $\gamma \in \DJs$ and $\beta \in \DJp$}.
\end{equation}
%

%
\subsection{Conjecture due to Kato and Sagaki.}
\label{subsec:conj}

Let $w,x \in W$, and $d \in Q^{\vee,+}$.
Following \cite{LNSX} and \cite{KS}, we define $\bR{w,x,d}$ to be the set of all directed paths
\begin{equation} \label{eq:dpR}
\bp: w = u_{0} \edge{\beta_{1}} u_{1} \edge{\beta_{2}} \cdots
\edge{\beta_{r}} u_{r} = \ed(\bp)
\end{equation}
in $\QBG(W)$ starting at $w$ satisfying the conditions that
%
%
\begin{equation} \label{eq:R}
\begin{cases}
r \ge 0, \quad
\beta_{1} \lhd \cdots \lhd \beta_{r}, \quad
\beta_{t} \in \Delta^{+}_{2} \ \text{for all $1 \le t \le r$}, \\[1mm]
\qwt (\ed(\bp) \Rightarrow x) \le d - \qwt(\bp), \quad
\pair{\vpi_{i}}{d - \qwt(\bp)} = 0, \quad \ed(\bp) \in x\WJs.
\end{cases}
\end{equation}
Note that an element $\bp \in \bR{w,x,d}$ is a label-increasing
(and hence shortest, by \cite[Theorem~7.3(2)]{LNSSS1}) directed path from $w$ to $\ed(\bp)$.
Thus,
\begin{equation} \label{eq:Rshort}
\ell(\bp) = \ell(w \Rightarrow \ed(\bp)) \qquad \text{and} \qquad
\qwt(\bp) = \qwt(w \Rightarrow \ed(\bp)).
\end{equation}

Here we recall the following conjecture due to Kato and Sagaki.
%
%
\begin{conj}[{\cite[Conjecture~3.5]{KS}}] \label{conj:R1even}
Assume that $\vpi_{i}$ is quasi-minuscule.
Then we have $\# \bR{w,x,d} \in \{ 1 \} \sqcup 2\BZ_{\ge 0}$
for all $w,x \in W$ and $d \in Q^{\vee,+}$.
\end{conj}

The following is the main result of this paper.
%
%
\begin{thm} \label{thm:main}
Conjecture~\ref{conj:R1even} is true.
Moreover, if $\# \bR{w,x,d}$ is even, then
\begin{equation} \label{eq:main2}
\sum_{ \bp \in \bR{w,x,d} } (-1)^{\ell(\bp)} = 0.
\end{equation}
\end{thm}

The proof of Theorem~\ref{thm:main} will be completed in Section~\ref{sec:prfmain}.
\begin{rem} \label{rem:main3}
It is proved in \cite[Theorem~3.7]{KS} that if Conjecture~\ref{conj:R1even} is true, then
\eqref{eq:main2} holds whenever $\# \bR{w,x,d} \ge 2$.
Our proof of \eqref{eq:main2} is independent of \cite[Theorem~3.7]{KS};
we will obtain $\# \bR{w,x,d} \in \{ 1 \} \sqcup 2\BZ_{\ge 0}$ and \eqref{eq:main2} simultaneously.
\end{rem}

%
\subsection{Cancellation-free version of the quantum $K$-theoretic divisor axiom.}
\label{subsec:cfree}

Fix a Borel subgroup $B$ such that $T \subset B \subset G$,
and let $X := G/B$ be the flag manifold.
For $u \in W$, let $\CO_u$ and $\CO^u$ denote the Schubert class
and the opposite Schubert class in the $T$-equivariant $K$-theory ring $K_T(X)$,
respectively.
For $\gamma_1,\ldots,\gamma_m \in K_T(X)$ and $d \in Q^{\vee,+}$,
we denote by $\langle\gamma_1,\ldots,\gamma_m\rangle_d$
the corresponding $m$-point ($T$-equivariant) $K$-theoretic Gromov-Witten
(KGW) invariant, which takes values in $K_T(\pt) \cong \BZ[\Lambda]$;
see \cite[Section~2.4]{KS} and \cite{LNSX} for details.
We do not use the geometric definition of these invariants.
%
%
\begin{prop}[{\cite[Lemma~4.1]{LNSX}}] \label{prop:twp}
Let $d \in Q^{\vee,+}$, and $w,x \in W$.
Then, we have
\begin{equation} \label{eq:2pt-qwt}
\langle \CO^{w}, \CO_{x} \rangle_{d} =
\begin{cases}
1 & \text{\rm if $\qwt(w \Rightarrow x) \le d$}, \\
0 & \text{\rm otherwise}.
\end{cases}
\end{equation}
\end{prop}

Let $i \in I$ be such that $\vpi_i$ is quasi-minuscule; recall that $\J = I \setminus \{i\}$.
We know from \cite[Theorem~3.2]{LNSX} and \cite[(3.4)]{KS} that
%
%
\begin{equation} \label{eq:thp-qm}
\langle \CO^{s_i}, \CO^{w}, \CO_{x} \rangle_{d} =
 \langle \CO^{w}, \CO_{x} \rangle_{d} - \be^{-\vpi_{i}+\wt(\eta)}
 \sum_{ \bp \in \bR{w,x,d} }
 (-1)^{\ell(\bp)},
\end{equation}
where $\wt(\eta) := \frac{1}{2}w\vpi_{i}+\frac{1}{2}x\vpi_{i} \in \Lambda$
when $\bR{w,x,d}\ne\emptyset$
(see \cite[Section~3.1]{LNSX} and \cite[Section~3.2]{KS} for the notation $\eta$).
%
%
Combining Theorem~\ref{thm:main} with \cite[Theorems~3.4 and 3.7]{KS},
we obtain the following corollary.
\begin{cor} \label{cor:cfree}
Assume that $\vpi_{i}$ is quasi-minuscule; recall that $J = I \backslash \{ i \}$.
Let $w,x \in W$, and $d \in Q^{\vee,+}$.
Setting $x_{\min}:=\tbmin{x}{\J}{w}$, we have
\begin{equation} \label{eq:cfree}
\langle \CO^{s_i}, \CO^{w}, \CO_{x} \rangle_{d} =
\begin{cases}
1 - (-1)^{\ell(w) - \ell(x_{\min})} \be^{-\vpi_{i}+\wt(\eta)}
  & \text{\rm if $\# \bR{w,x,d} = 1$}, \\[1mm]
1 & \text{\rm if $\# \bR{w,x,d} \ge 2$}, \\[1mm]
1 & \text{\rm if $\bR{w,x,d} = \emptyset$ and $\qwt(w \Rightarrow x) \le d$}, \\[1mm]
0 & \text{\rm if $\bR{w,x,d} = \emptyset$ and $\qwt(w \Rightarrow x) \not\le d$}.
\end{cases}
\end{equation}
\end{cor}

%
\subsection{Reduction of Theorem~\ref{thm:main}.}
\label{subsec:E}

In order to prove Theorem~\ref{thm:main},
we describe the set of endpoints $\ed(\bp)$ of directed paths in $\bR{w,x,d}$
in terms of the dual tilted Bruhat order (see Section~\ref{subsec:tilted}).
For an arbitrary subset $\J$ of $I$, and for $w,x \in W$
and $d \in Q^{\vee,+}$, we set
%
%
\begin{equation} \label{eq:E}
\bE{w,x,d} = \bE{w,x,d}^{\J} :=
\bigl\{ v \in x\WJs \mid \tbmax{w}{\J}{v} = w, \
 \qwt(w \Rightarrow v) + \qwt(v \Rightarrow x) \le d \bigr\}.
\end{equation}
\begin{prop} \label{prop:end}
Let $i \in I$ be such that $\vpi_{i}$ is quasi-minuscule,
and set $\J = I \setminus \{i\}$.
Let $w,x \in W$, and $d \in Q^{\vee,+}$.
If $\bR{w,x,d} \ne \emptyset$, then the map
\begin{equation} \label{eq:end}
\ed : \bR{w,x,d} \rightarrow \bE{w,x,d}, \qquad \bp \mapsto \ed(\bp),
\end{equation}
is bijective.
\end{prop}

\begin{proof}
First, let us show that $\ed(\bp) \in \bE{w,x,d}$ for all $\bp \in \bR{w,x,d}$.
Let $\bp \in \bR{w,x,d}$, and set $v := \ed(\bp) \in x\WJs$.
Because the labels of $\bp$ are all contained in $\Delta^{+}_{2} \subset \DJp$,
it follows from \cite[Lemma~2.15]{KNS} that $\tbmax{w}{\J}{v} = w$.
Also, by \eqref{eq:Rshort} and \eqref{eq:R}, we have
$\qwt(w \Rightarrow v) + \qwt(v \Rightarrow x) = \qwt(\bp) + \qwt(\ed(\bp) \Rightarrow x) \le d$.
Hence we obtain $v \in \bE{w,x,d}$.
The map \eqref{eq:end} is injective, since a label-increasing
directed path from $w$ is uniquely determined by its endpoint
\cite[Theorem~7.3(1)]{LNSSS1}.

Next, let us show that the map \eqref{eq:end} is surjective.
Since $\bR{w,x,d} \ne \emptyset$ by the assumption,
we can take (and fix) an arbitrary element $\bp_{0} \in \bR{w,x,d}$; we set $v_{0} := \ed(\bp_{0}) \in x\WJs$.
Let $v \in \bE{w,x,d}$.
Let $\bq$ be the label-increasing directed path from $w$ to $v$
in $\QBG(W)$.
We show that $\bq \in \bR{w,x,d}$.

We first claim that
\begin{equation} \label{eq:end1}
\text{the labels of $\bq$ are all contained in $\Delta^{+}_{2}$}.
\end{equation}
Indeed, since $v_{0},v \in x\WJs$, it follows from
\cite[Lemma~2.9]{KS} that there exists a directed path
$v_{0} \edge{\gamma_{1},\dots,\gamma_{t}} v$
with $\gamma_{1},\dots,\gamma_{t} \in \DJs$.
Then the concatenation
\begin{equation*}
\underbrace{w \edge{\beta_{1},\dots,\beta_{r}} v_{0}}_{=\bp_{0}}
\edge{\gamma_{1},\dots,\gamma_{t}} v
\end{equation*}
is a directed path from $w$ to $v$ lying in
$\QBG_{(1/2)\vpi_{i}}(W)$.
By \cite[Theorem~7.3]{LNSSS1}, $\bq$ is a shortest directed path
from $w$ to $v$. Hence, by \cite[Lemma~6.7]{LNSSS2},
$\bq$ also lies in $\QBG_{(1/2)\vpi_i}(W)$.
Because $\vpi_{i}$ is quasi-minuscule, we obtain
$\pair{\vpi_{i}}{\beta^{\vee}} \in \{0,2\}$
for every label $\beta$ in $\bq$.
Also, since $\tbmax{w}{\J}{v} = w$ by the definition of
$\bE{w,x,d}$, it follows from \cite[Lemma~2.15]{KNS}
that all the labels in $\bq$ belong to
$
\DJp =
\bigl\{ \beta \in \Delta^{+}
\mid \pair{\vpi_{i}}{\beta^{\vee}} \ne 0 \bigr\}.
$
Thus we have shown \eqref{eq:end1}.

Since $\bq$ is a shortest directed path from $w$ to $v$,
we have $\qwt(\bq) = \qwt(w \Rightarrow v)$.
Hence, by the definition of $\bE{w,x,d}$, we obtain
\begin{equation*}
\qwt(\bq) + \qwt(\ed(\bq) \Rightarrow x)
= \qwt(w \Rightarrow v) + \qwt(v \Rightarrow x)
\le d.
\end{equation*}
Moreover, since $v_{0},v \in x\WJs$, it follows from
Lemma~\ref{lem:const} that
\begin{equation*}
\begin{aligned}
\pair{\vpi_{i}}{\qwt(\bq)}
&= \pair{\vpi_{i}}{\qwt(w \Rightarrow v)}
 = \pair{\vpi_{i}}{\qwt(w \Rightarrow v_{0})} \\
&= \pair{\vpi_{i}}{\qwt(\bp_{0})}
= \pair{\vpi_{i}}{d}.
\end{aligned}
\end{equation*}
Therefore we conclude that $\bq \in \bR{w,x,d}$.
Since $\ed(\bq) = v$, the map \eqref{eq:end} is surjective.
\end{proof}

By Proposition~\ref{prop:end}, in order to prove Theorem~\ref{thm:main}
(which is the main result of this paper)
it suffices to prove that
for an arbitrary subset $\J$ of $I$ and $w,x \in W$, $d \in Q^{\vee,+}$ satisfying
\begin{equation} \label{eq:asm}
\pair{\vpi_{k}}{d} = \pair{\vpi_{k}}{\qwt(w \Rightarrow x)}
\qquad \text{for all $k \in I \setminus \J$},
\end{equation}
the following holds:
\begin{equation} \label{eq:Eclaim}
\text{if $\# \bE{w,x,d} \ge 2$, then} \quad
\sum_{v \in \bE{w,x,d}} (-1)^{\ell(v)} = 0.
\end{equation}
Indeed, we apply the claim above to the case
where $\J = I \setminus \{i\}$ with $i \in I$ such that $\varpi_i$ is quasi-minuscule
and $\# \bR{w,x,d} \ge 2$ with $w,x \in W$, $d \in Q^{\vee,+}$.
For $\bp \in \bR{w,x,d}$, we see by \eqref{eq:R}, \eqref{eq:Rshort},
and Lemma~\ref{lem:const} give that
$
\pair{\vpi_i}{d}
= \pair{\vpi_i}{\qwt(\bp)}
= \pair{\vpi_i}{\qwt(w \Rightarrow x)},
$
since $\ed(\bp) \in x\WJs$.
Hence \eqref{eq:asm} holds.
By Proposition~\ref{prop:end} and \eqref{eq:Eclaim}, we obtain
\begin{equation*}
\sum_{\bp \in \bR{w,x,d}} (-1)^{\ell(\bp)}
=
(-1)^{\ell(w)}
\sum_{v \in \bE{w,x,d}} (-1)^{\ell(v)}
= 0,
\end{equation*}
which implies that $\# \bR{w,x,d}$ is an even number.
Thus we have proved Theorem~\ref{thm:main}.
We will prove \eqref{eq:Eclaim} under assumption \eqref{eq:asm}
in Theorem~\ref{thm:E}
in terms of the Bruhat order and the Demazure products
on the affine Weyl group.

%
\section{Bruhat order and Demazure products on the affine Weyl group.}
\label{sec:aff}

In this section, we fix our notation for the affine Weyl group $\Waf$,
and recall from \cite{Sch} the descriptions of the Bruhat order and
of Demazure products on $\Waf$ in terms of the quantum Bruhat graph $\QBG(W)$.
We then prove Propositions~\ref{prop:Bruhat} and \ref{prop:trans},
which will be used in Section~\ref{sec:proof}.

%
\subsection{Affine Weyl group.}
\label{subsec:affine}

We recall the affine root system and the affine Weyl group; for details,see, e.g., \cite[Section~6]{Kac}.

Let
\begin{equation*}
\rr := \bigl\{ \alpha + n\delta \mid \alpha \in \Delta, \, n \in \BZ \bigr\}, \qquad
\prr := \Delta^{+} \sqcup \bigl\{ \alpha + n\delta \mid \alpha \in \Delta, \, n \in \BZ_{> 0} \bigr\}
\end{equation*}
be the set of real roots and the set of positive real roots
of the untwisted affine root system associated to $\Delta$, respectively.
Here $\delta$ is the null root.
We set $\alpha_{0} := -\theta + \delta \in \prr$, and $\Iaf := I \sqcup \{0\}$.
For $\xi \in Q^{\vee}$, let $t_{\xi}$ denote the translation by $\xi$.
The (affine) Weyl group $\Waf$ is the semi-direct product
$\Waf = W \ltimes \bigl\{ t_{\xi} \mid \xi \in Q^{\vee} \bigr\}$
of the (finite) Weyl group $W$ and the abelian group $\bigl\{ t_{\xi} \mid \xi \in Q^{\vee} \bigr\} \cong Q^{\vee}$,
in which $v t_{\xi} v^{-1} = t_{v\xi}$ for $v \in W$ and $\xi \in Q^{\vee}$,
and $t_{\xi}t_{\zeta} = t_{\xi+\zeta}$ for $\xi ,\zeta \in Q^{\vee}$. 
Hence we have
%
%
\begin{equation} \label{eq:prod}
(vt_{\xi})(ut_{\zeta}) = vu\,t_{u^{-1}\xi+\zeta} \qquad
\text{for $v,u \in W$ and $\xi,\zeta \in Q^{\vee}$}.
\end{equation}
The group $\Waf$ acts on $\rr$ by
%
%
\begin{equation} \label{eq:act}
vt_{\xi}(\alpha+n\delta) = v\alpha + \bigl( n-\pair{\alpha}{\xi} \bigr) \delta \qquad
\text{for $v \in W$, $\xi \in Q^{\vee}$, and $\alpha+n\delta \in \rr$}.
\end{equation}
For $\beta = \alpha+n\delta \in \rr$, we set
$s_{\beta} := s_{\alpha}t_{n\alpha^{\vee}} \in \Waf$;
in particular, $s_{0} := s_{\alpha_{0}} = s_{\theta}t_{-\theta^{\vee}}$.
Then, $\Waf$ is a Coxeter group with $\{ s_{j} \mid j \in \Iaf \}$
the set of simple reflections.
Denote by $e$ the identity element of $\Waf$.
We also denote by $\ell : \Waf \rightarrow \BZ_{\ge 0}$
the length function on $\Waf$, and by $\le$ the Bruhat order on $\Waf$;
since $W = \langle s_{j} \mid j \in I \rangle$ is a (standard) parabolic subgroup of $\Waf$,
their restrictions to $W$ agree with on $W$.
Recall that
\begin{equation} \label{eq:lenaf}
\ell(y) = \# \bigl\{ \beta \in \prr \mid y\beta \in -\prr \bigr\} \qquad
\text{for $y \in \Waf$}.
\end{equation}
We use the Demazure product $\ast$ on the Coxeter group $\Waf$
defined in Section~\ref{subsec:Demazure}.

%
\subsection{Length positive elements.}
\label{subsec:LP}

For $\beta \in \Delta$, we set $\chi(\beta) := 1$ if $\beta \in \Delta^{+}$, and
$\chi(\beta) := 0$ if $\beta \in \Delta^{-}$.
%
%
\begin{dfn}[{\cite[Definition~2.2]{Sch}}] \label{dfn:LP}
Let $y = vt_{\xi} \in \Waf$, with $v \in W$ and $\xi \in Q^{\vee}$.
\begin{enu}
\item For $\alpha \in \Delta$, we set
$\ell(y,\alpha) := \pair{\alpha}{\xi} + \chi(\alpha) - \chi(v\alpha) \in \BZ$.

\item An element $u \in W$ is said to be length positive for $y$ if
$\ell(y,u\alpha) \ge 0$ for all $\alpha \in \Delta^{+}$.
We denote by $\LP(y)$ the set of all length positive elements for $y$,
that is,
\[
\LP(y)
:=
\bigl\{ u \in W \mid
\ell(y,u\alpha) \ge 0
\text{ for all } \alpha \in \Delta^{+} \bigr\}.
\]
\end{enu}
\end{dfn}

Since $\chi(-\alpha) = 1 - \chi(\alpha)$ for $\alpha \in \Delta$, we have
\begin{equation} \label{eq:lfun}
\ell(y,-\alpha) = -\ell(y,\alpha) \qquad
\text{for $y \in \Waf$ and $\alpha \in \Delta$}.
\end{equation}

An element $\xi \in Q^{\vee}$ is said to be dominant (resp., regular dominant)
if $\pair{\alpha_{j}}{\xi} \ge 0$ (resp., $\pair{\alpha_{j}}{\xi} > 0$) for all $j \in I$.

%
%
\begin{lem} \label{lem:LP}
\mbox{}
\begin{enu}
\item Let $J$ be an arbitrary subset of $I$,
and let $\xi \in Q^{\vee}$ be a dominant element,
such that
$\J = \bigl\{ j \in I \mid \pair{\alpha_{j}}{\xi} = 0 \bigr\}$.
Then, $\LP(t_{v\xi}) = v\WJs$ for all $v \in W$.

\item Let $\zeta \in Q^{\vee}$ be a dominant element. Then, $e \in \LP(vt_{\zeta})$ for all $v \in W$.

\item Let $\zeta \in Q^{\vee}$ be a regular dominant element. Then,
$\LP(vt_{\zeta}) = \{e\}$ and $\ell(vt_{\zeta}) = \ell(v) + 2\pair{\rho}{\zeta}$ for all $v \in W$.
\end{enu}
\end{lem}

\begin{proof}
(1) Let $u \in W$. Since $\ell(t_{v\xi},u\alpha) = \pair{u\alpha}{v\xi} = \pair{\alpha}{u^{-1}v\xi}$
for $\alpha \in \Delta^{+}$, we see that $u \in \LP(t_{v\xi})$ if and only if
$u^{-1}v\xi$ is dominant.
Because $\xi$ is the unique dominant element in the $W$-orbit of $\xi$, and
the stabilizer of $\xi$ in $W$ is identical to $\WJs$,
we deduce that $u^{-1}v\xi$ is dominant if and only if $u^{-1}v \in \WJs$,
or equivalently, $u \in v\WJs$.

(2) For $\alpha \in \Delta^{+}$, we have
$\ell(vt_{\zeta},\alpha) = \pair{\alpha}{\zeta} + 1 - \chi(v\alpha) \ge \pair{\alpha}{\zeta} \ge 0$.
Hence, $e \in \LP(vt_{\zeta})$.

(3) Since $\zeta$ is regular dominant, we see by the computation in the proof of part (2) that
$\ell(vt_{\zeta},\alpha) > 0$ for all $\alpha \in \Delta^{+}$, and hence
 by \eqref{eq:lfun}, $\ell(vt_{\zeta},\alpha) < 0$ for all $\alpha \in \Delta^{-}$.
If $u \in W$ and $u \ne e$, then there exists $\alpha \in \Delta^{+}$ such that $u\alpha \in \Delta^{-}$,
which implies that $u \notin \LP(vt_{\zeta})$. 
Thus $\LP(vt_{\zeta}) = \{e\}$.
It follows from
\cite[Lemma~2.3]{Sch} that
$\ell(vt_{\zeta}) = \ell(v) + 2\pair{\rho}{\zeta}$.
\end{proof}

%
\subsection{Results by Schremmer and their consequences.}
\label{subsec:Sch}

We recall two results by Schremmer \cite{Sch}.
The first one is a criterion for the Bruhat order on $\Waf$.
%
%
\begin{thm}[{\cite[Theorem~4.2]{Sch}}] \label{thm:Sch1}
Let $y_{1} = v_{1}t_{\xi_{1}}$ and $y_{2} = v_{2}t_{\xi_{2}}$
be elements of $\Waf$,
with $v_{1},v_{2} \in W$ and $\xi_{1},\xi_{2} \in Q^{\vee}$,
and let $u_{1} \in \LP(y_{1})$.
Then, $y_{1} \le y_{2}$ in the Bruhat order on $\Waf$
if and only if there exists $u_{2} \in W$ such that
\begin{equation} \label{eq:Sch1}
u_{1}^{-1}\xi_{1}
+ \qwt(u_{2} \Rightarrow u_{1})
+ \qwt(v_{1}u_{1} \Rightarrow v_{2}u_{2})
\le u_{2}^{-1}\xi_{2}.
\end{equation}
\end{thm}
\begin{rem}
In \cite[Theorem~4.2]{Sch}, a criterion for the Bruhat order
on (extended) affine Weyl groups is given in terms of
a ``Bruhat-deciding datum'' $(u_{1},J_{1},\dots,J_{m})$
for $y_{1}$ (see \cite[Definition~4.1]{Sch}).
Theorem~\ref{thm:Sch1} above is the special case of this criterion
in which $m = 1$ and $J_{1} = \emptyset$;
note that $(u_{1},\emptyset)$ is a Bruhat-deciding datum
for $y_{1}$ if (and only if) $u_{1} \in \LP(y_{1})$.
\end{rem}

The second one is an explicit formula for the Demazure product on $\Waf$.
%
%
\begin{thm}[{\cite[Theorem~5.11]{Sch}}] \label{thm:Sch2}
Let $y_{1} = v_{1}t_{\xi_{1}}$ and $y_{2} = v_{2}t_{\xi_{2}}$ be elements of $\Waf$,
with $v_{1},v_{2} \in W$ and $\xi_{1},\xi_{2} \in Q^{\vee}$.
Let $(u_{1},u_{2}) \in \LP(y_{1}) \times \LP(y_{2})$ be such that
\begin{equation*}
\ell(u_{1} \Rightarrow v_{2}u_{2}) = \min \bigl\{ \ell(u_{1}' \Rightarrow v_{2}u_{2}') \mid
 (u_{1}',u_{2}') \in \LP(y_{1}) \times \LP(y_{2}) \bigr\}.
\end{equation*}
Then the Demazure product $y_{1} \ast y_{2} \in \Waf$ is given by
\begin{equation} \label{eq:Sch2}
y_{1} \ast y_{2} =
v_{1}u_{1}\,t_{u_{1}^{-1}\xi_{1}+u_{2}^{-1}\xi_{2}-\qwt(u_{1} \Rightarrow v_{2}u_{2})}\,u_{2}^{-1}.
\end{equation}
\end{thm}

\begin{prop} \label{prop:Bruhat}
Let $v_{1},v_{2} \in W$.
Let $\zeta_{1} \in Q^{\vee}$ be a dominant element, and
let $\zeta_{2} \in Q^{\vee}$ be a regular dominant element. Then,
\begin{equation*}
v_{1}t_{\zeta_{1}} \le v_{2}t_{\zeta_{2}}
\iff \zeta_{1} + \qwt(v_{1} \Rightarrow v_{2}) \le \zeta_{2}.
\end{equation*}
\end{prop}
\begin{proof}
Since $e \in \LP(v_{1} t_{\zeta_{1}})$
by Lemma~\ref{lem:LP}\,(2),
it follows from Theorem~\ref{thm:Sch1} that
$v_{1} t_{\zeta_{1}} \le v_{2} t_{\zeta_{2}}$ if and only if
there exists $u \in W$ such that
\begin{equation} \label{eq:B1}
\zeta_{1} + \qwt(u \Rightarrow e)
+ \qwt(v_{1} \Rightarrow v_{2} u)
\le u^{-1}\zeta_{2}.
\end{equation}
If $\zeta_{1} + \qwt(v_{1} \Rightarrow v_{2}) \le \zeta_{2}$,
then \eqref{eq:B1} holds for $u = e$.
Conversely, assume that \eqref{eq:B1} holds for some $u \in W$;
we show by descending induction on $\ell(u)$ that \eqref{eq:B1}
holds also for $u = e$.
If $u = e$, then there is nothing to prove.
Assume that $u \ne e$, and take $j \in I$ such that $us_{j} < u$.
We have
\begin{align*}
& \qwt(us_{j} \Rightarrow e)
\le \qwt(us_{j} \Rightarrow u) + \qwt(u \Rightarrow e)
= \qwt(u \Rightarrow e), \\
& \qwt(v_{1} \Rightarrow v_{2}us_{j})
\le \qwt(v_{1} \Rightarrow v_{2}u)
+ \qwt(v_{2}u \Rightarrow v_{2}us_j)
\le \qwt(v_{1} \Rightarrow v_{2}u) + \alpha_{j}^{\vee}.
\end{align*}
Also, since $us_{j} < u$, we have $u\alpha_{j} \in \Delta^{-}$,
and hence $\pair{u\alpha_{j}}{\zeta_{2}} \le -1$
since $\zeta_{2}$ is regular dominant. Therefore,
\begin{equation*}
(us_{j})^{-1}\zeta_{2} = s_{j}u^{-1}\zeta_{2}
= u^{-1}\zeta_{2}
- \pair{u\alpha_{j}}{\zeta_{2}}\alpha_{j}^{\vee}
\ge u^{-1}\zeta_{2} + \alpha_{j}^{\vee}.
\end{equation*}
Combining these inequalities with \eqref{eq:B1}, we obtain
\begin{align*}
\zeta_{1} + \qwt(us_{j} \Rightarrow e)
+ \qwt(v_{1} \Rightarrow v_{2}us_{j})
& \le \zeta_{1} + \qwt(u \Rightarrow e)
+ \qwt(v_{1} \Rightarrow v_{2}u) + \alpha_{j}^{\vee} \\
& \le u^{-1}\zeta_{2} + \alpha_{j}^{\vee}
\le (us_{j})^{-1}\zeta_{2},
\end{align*}
as desired. This proves the proposition.
\end{proof}
%
%
%
\begin{prop} \label{prop:trans}
Let $\xi \in Q^{\vee}$ be a dominant element, and set
$\J := \bigl\{ j \in I \mid \pair{\alpha_{j}}{\xi} = 0 \bigr\}$.
Let $v,u \in W$, and let $\zeta \in Q^{\vee}$
be a regular dominant element. Then,
\begin{equation} \label{eq:trans}
t_{v\xi} \ast (ut_{\zeta})
= u_{1}t_{\xi+\zeta-\qwt(u_{1} \Rightarrow u)}, \qquad
\text{\rm where $u_{1} := \tbmax{v}{\J}{u}$}.
\end{equation}
\end{prop}
\begin{proof}
We know from Lemma~\ref{lem:LP}\,(1) and (3) that
$\LP(t_{v\xi}) = v\WJs$ and $\LP(ut_{\zeta}) = \{e\}$.
Also, it follows from \eqref{eq:tbmax1} that
\begin{equation*}
\ell(u_{1} \Rightarrow u)
= \min \bigl\{ \ell(v' \Rightarrow u)
\mid v' \in v\WJs \bigr\}.
\end{equation*}
Therefore, by applying Theorem~\ref{thm:Sch2} to
$y_{1} = t_{v\xi}$, $y_{2} = ut_{\zeta}$
and the pair $(u_{1},e)$, we obtain
\begin{equation*}
t_{v\xi} \ast (ut_{\zeta})
= u_{1}t_{u_{1}^{-1}v\xi+\zeta-\qwt(u_{1} \Rightarrow u)}.
\end{equation*}
Since $u_{1}^{-1}v \in \WJs$,
and since $\WJs$ stabilizes $\xi$,
we have $u_{1}^{-1}v\xi = \xi$.
This proves the proposition.
\end{proof}

%
\section{Bijection with subset of affine Weyl group.}
\label{sec:proof}

Throughout this section, we fix an arbitrary subset $\J$ of $I$,
and let $w,x \in W$ and $d \in Q^{\vee,+}$ satisfy \eqref{eq:asm}.
Recall the set $\bE{w,x,d} \subset x\WJs$ from \eqref{eq:E}, and
the sets $\CE_u(f,g)$ from \eqref{eq:CE}.
In this section, we apply the definition of $\CE_u(f,g)$
given in Section~2 to the affine Weyl group $W_{\mathrm{af}}$.

We fix a dominant element $\xi \in Q^{\vee}$ such that
$\bigl\{ j \in I \mid \pair{\alpha_{j}}{\xi} = 0 \bigr\} = \J$.
We define an element $z \in \Waf$ by
\begin{equation} \label{eq:x0}
z := t_{w\xi} \ast t_{x\xi}.
\end{equation}
\begin{lem} \label{lem:x0}
Let $v \in W$, and set
\begin{equation} \label{eq:hvhw}
\ha{v} := \tbmax{x}{\J}{v} \in x\WJs, \qquad
\ha{w} := \tbmax{w}{\J}{\ha{v}} \in w\WJs.
\end{equation}
Let $\zeta \in Q^{\vee}$ be such that both $\zeta$ and
$\zeta - \qwt(\ha{v} \Rightarrow v)$ are regular dominant. Then,
\begin{equation} \label{eq:x0v}
z \ast (vt_{\zeta}) =
\ha{w}\,t_{2\xi+\zeta-\qwt(\ha{v} \Rightarrow v)-\qwt(\ha{w} \Rightarrow \ha{v})}.
\end{equation}
\end{lem}

\begin{proof}
By the associativity of the Demazure product $\ast$, we have
$z \ast (vt_{\zeta}) =
t_{w\xi} \ast \bigl( t_{x\xi} \ast (vt_{\zeta}) \bigr)$.
Since $\zeta$ is regular dominant, it follows from
Proposition~\ref{prop:trans} that
$t_{x\xi} \ast (vt_{\zeta}) = \ha{v}\,t_{\zeta_{1}}$, where
$\zeta_{1} := \xi + \zeta - \qwt(\ha{v} \Rightarrow v)$.
Since $\xi$ is dominant and
$\zeta - \qwt(\ha{v} \Rightarrow v)$ is regular dominant,
we see that $\zeta_{1}$ is regular dominant.
Applying Proposition~\ref{prop:trans} again, we obtain
$t_{w\xi} \ast (\ha{v}\,t_{\zeta_{1}})
= \ha{w}\,t_{\xi+\zeta_{1}-\qwt(\ha{w} \Rightarrow \ha{v})}$.
This proves the lemma.
\end{proof}

Now, let us fix $\eta \in Q^{\vee}$ as follows.
Since the subset $\{ a \in \Waf \mid a \le z \}$ is a finite set,
the subset
%
%
\begin{equation} \label{eq:Xi}
\begin{split}
\Xi :=
& \bigl\{ \xi_{a} - \qwt(v' \Rightarrow v'') \mid a \in \Waf, \ a \le z, \ v',v'' \in W \bigr\} \\
& \cup \bigl\{ \qwt(v' \Rightarrow v'') \mid v',v'' \in W \bigr\} \cup \{d\}
\end{split}
\end{equation}
of $Q^{\vee}$ is also a finite set, where for $a \in \Waf$, we define
$v_{a} \in W$ and $\xi_{a} \in Q^{\vee}$ by $a^{-1}wt_{2\xi} = v_{a}t_{\xi_{a}}$;
note that $\xi_{a} = \xi_{a} - \qwt(e \Rightarrow e) \in \Xi$
for all $a \in \Waf$ such that $a \le z$.
We fix $\eta \in Q^{\vee}$ such that
%
%
\begin{equation} \label{eq:eta}
\text{$\eta + \gamma$ is regular dominant for all $\gamma \in \Xi$}.
\end{equation}
Then we define $\ti{w},\ti{x} \in \Waf$ by
%
%
\begin{equation} \label{eq:zu}
\ti{w} := wt_{2\xi+\eta}, \qquad \ti{x} := xt_{\eta+d}.
\end{equation}
%
%
\begin{prop} \label{prop:lift}
Keep the notation and setting above; in particular, we assume \eqref{eq:asm}.
Recall that
$\CE_{z}(\ti{w},\ti{x})
= \bigl\{ y \in \Waf \mid y \le \ti{x},\ z \ast y = \ti{w} \bigr\}$.
The map
\begin{equation} \label{eq:Phi}
\Phi \colon \bE{w,x,d} \longrightarrow \CE_{z}(\ti{w},\ti{x}),
\qquad
v \longmapsto vt_{\eta+\qwt(w \Rightarrow v)}
\end{equation}
is bijective and satisfies
$(-1)^{\ell(\Phi(v))} = (-1)^{\ell(v)}$
for all $v \in \bE{w,x,d}$.
\end{prop}

\begin{proof}
First, let us show that $\Phi(v) := vt_{\eta+\qwt(w \Rightarrow v)} \in \CE_{z}(\ti{w},\ti{x})$
for all $v \in \bE{w,x,d}$.
Let $v \in \bE{w,x,d}$, and set $\zeta := \eta + \qwt(w \Rightarrow v)$;
note that $\zeta$ is regular dominant by \eqref{eq:eta} since $\qwt(w \Rightarrow v) \in \Xi$.
Since $v \in x\WJs$, we have $\tbmax{x}{\J}{v} = v$.
Also, we have $\tbmax{w}{\J}{v} = w$ by the definition of $\bE{w,x,d}$.
By Lemma~\ref{lem:x0}, we have
\begin{equation*}
z \ast \Phi(v) = wt_{2\xi+\zeta-\qwt(w \Rightarrow v)} = wt_{2\xi+\eta} = \ti{w}.
\end{equation*}
Also, since $\zeta$ and $\eta+d$ are regular dominant by \eqref{eq:eta}
(recall that $d \in \Xi$), it follows from Proposition~\ref{prop:Bruhat} that
\begin{equation*}
\Phi(v) = vt_{\zeta} \le xt_{\eta+d} = \ti{x} \iff
\eta + \qwt(w \Rightarrow v) + \qwt(v \Rightarrow x) \le \eta + d.
\end{equation*}
The inequality on the right-hand side holds since $v \in \bE{w,x,d}$.
Hence we get $\Phi(v) \le \ti{x}$. Thus we have shown that $\Phi(v) \in \CE_{z}(\ti{w},\ti{x})$.

It is obvious that the map $\Phi$ is injective. Let us show that $\Phi$ is surjective.
Let $y \in \CE_{z}(\ti{w},\ti{x})$, and set $a := \ti{w}y^{-1} \in \Waf$.
Since $z \ast y = \ti{w}$, we see by Proposition~\ref{prop:He} that $a \le z$.
Hence it follows that
$
y = a^{-1}\ti{w} = (a^{-1}wt_{2\xi})t_{\eta} = v_{a}t_{\xi_{a}+\eta}
$
in the notation of \eqref{eq:Xi}.
We set $v := v_{a} \in W$ and $\zeta := \xi_{a}+\eta \in Q^{\vee}$,
and define $\ha{v},\ha{w} \in W$ by
\begin{equation*}
\ha{v} := \tbmax{x}{\J}{v},
\qquad
\ha{w} := \tbmax{w}{\J}{\ha{v}}.
\end{equation*}
Since $\xi_{a} , \xi_{a} - \qwt(\ha{v} \Rightarrow v) \in \Xi$,
both $\zeta = \eta + \xi_{a}$ and
$\zeta - \qwt(\ha{v} \Rightarrow v) = \eta + \xi_{a} - \qwt(\ha{v} \Rightarrow v)$
are regular dominant by \eqref{eq:eta}.
By Lemma~\ref{lem:x0}, we have
\begin{equation*}
wt_{2\xi+\eta} = \ti{w} = z \ast y =
\ha{w}\,t_{2\xi+\zeta-\qwt(\ha{v} \Rightarrow v)-\qwt(\ha{w} \Rightarrow \ha{v})}.
\end{equation*}
Therefore, we obtain
\begin{equation} \label{eq:lift1}
\ha{w} = w \qquad \text{and} \qquad
\zeta = \eta + \qwt(w \Rightarrow \ha{v}) + \qwt(\ha{v} \Rightarrow v).
\end{equation}
Also, since $y = vt_{\zeta} \le \ti{x} = xt_{\eta+d}$, and since $\zeta$ and $\eta+d$ are regular dominant,
it follows from Proposition~\ref{prop:Bruhat} that $\zeta + \qwt(v \Rightarrow x) \le \eta + d$.
Combining this inequality with \eqref{eq:lift1}, we get
\begin{equation} \label{eq:lift2}
\qwt(w \Rightarrow \ha{v}) + \qwt(\ha{v} \Rightarrow v) + \qwt(v \Rightarrow x) \le d.
\end{equation}
We claim that $v \in x\WJs$.
Since $\ha{v} \in x\WJs$, we see by Lemma~\ref{lem:const}
and the assumption \eqref{eq:asm} that
$\pair{\vpi_{k}}{\qwt(w \Rightarrow \ha{v})}
= \pair{\vpi_{k}}{\qwt(w \Rightarrow x)}
= \pair{\vpi_{k}}{d}$
for all $k \in I \setminus \J$.
Hence it follows from \eqref{eq:lift2} that
$\pair{\vpi_{k}}{\qwt(\ha{v} \Rightarrow v) + \qwt(v \Rightarrow x)} \le 0$.
Because $\qwt(\ha{v} \Rightarrow v) + \qwt(v \Rightarrow x) \in Q^{\vee,+}$, we deduce that
$\pair{\vpi_{k}}{\qwt(\ha{v} \Rightarrow v) + \qwt(v \Rightarrow x)} = 0$ for all $k \in I \setminus \J$,
that is, $\qwt(\ha{v} \Rightarrow v) + \qwt(v \Rightarrow x) \in \QJvp$.
Therefore, we obtain $v \in x\WJs$ by Lemma~\ref{lem:stay}, as desired.
Since $v \in x\WJs$, we have $\ha{v} = \tbmax{x}{\J}{v} = v$.
Then we see by \eqref{eq:lift1} and \eqref{eq:lift2} that
\begin{equation*}
\tbmax{w}{\J}{v} = \ha{w} = w, \qquad
\zeta = \eta + \qwt(w \Rightarrow v), \qquad
\qwt(w \Rightarrow v) + \qwt(v \Rightarrow x) \le d.
\end{equation*}
Hence $v \in \bE{w,x,d}$ and $y = \Phi(v)$.
This proves the surjectivity of $\Phi$.

Finally, let $v \in \bE{w,x,d}$. Since $\eta + \qwt(w \Rightarrow v)$ is regular dominant,
we see by Lemma~\ref{lem:LP}\,(3) that
$\ell(\Phi(v)) = \ell(v) + 2\pair{\rho}{\eta + \qwt(w \Rightarrow v)}$.
Since $\pair{\rho}{\gamma} \in \BZ$ for all $\gamma \in Q^{\vee}$,
we obtain $(-1)^{\ell(\Phi(v))} = (-1)^{\ell(v)}$.
This completes the proof of Proposition~\ref{prop:lift}.
\end{proof}

%
\section{Proof of Theorem~\ref{thm:main}.}
\label{sec:prfmain}

Let $i \in I$ be such that $\vpi_i$ is quasi-minuscule, and
set $\J := I \setminus \{i\}$ as in Section~\ref{sec:main}.
As shown at the end of Section~\ref{subsec:E},
Theorem~\ref{thm:main} follows from the assertion
\eqref{eq:Eclaim} with $J = I \backslash \{ i\}$.
Condition \eqref{eq:asm} holds whenever
$\bR{w,x,d} \ne \emptyset$, as shown there.

We now prove the following theorem for an arbitrary subset $\J$ of $I$
by combining Proposition~\ref{prop:lift} with Theorem~\ref{thm:cox}.%
%
\begin{thm} \label{thm:E}
Let $\J$ be an arbitrary subset of $I$.
Let $w,x \in W$ and $d \in Q^{\vee,+}$ satisfy \eqref{eq:asm}.
If $\# \bE{w,x,d} \ge 2$, then $\# \bE{w,x,d}$ is an even number, and
\begin{equation} \label{eq:Esum}
\sum_{v \in \bE{w,x,d}} (-1)^{\ell(v)} = 0.
\end{equation}
In particular, $\# \bE{w,x,d} \in \{1\} \sqcup 2\BZ_{\ge 0}$.
\end{thm}

\begin{proof}
Assume that $\# \bE{w,x,d} \ge 2$; in particular, $\bE{w,x,d} \ne \emptyset$.
Let $z,\ti{w},\ti{x} \in \Waf$ be as in \eqref{eq:x0} and \eqref{eq:zu}.
Proposition~\ref{prop:lift} gives that
\begin{equation*}
\# \bE{w,x,d} = \# \CE_{z}(\ti{w},\ti{x}) \qquad \text{and} \qquad
\sum_{v \in \bE{w,x,d}} (-1)^{\ell(v)} = \sum_{y \in \CE_{z}(\ti{w},\ti{x})} (-1)^{\ell(y)}.
\end{equation*}
Therefore, by applying Theorem~\ref{thm:cox} to
the Coxeter system $(\Waf,\,\{ s_{j} \mid j \in \Iaf \})$,
we deduce that $\# \bE{w,x,d}$ is an even number, and
$\sum_{v \in \bE{w,x,d}} (-1)^{\ell(v)} = 0$.
This completes the proof of Theorem~\ref{thm:E}.
\end{proof}

\begin{rem} \label{rem:single}
Keep the notation and setting of Proposition~\ref{prop:lift}.
We see by Theorem~\ref{thm:cox} and Proposition~\ref{prop:lift} that
$\# \bE{w,x,d} = 1$ if and only if $\ti{w} \dast \ti{x}^{-1} = z$;
in this case, the unique element $v$ of $\bE{w,x,d}$ is determined by
the equality $vt_{\eta+\qwt(w \Rightarrow v)} = z^{-1}\ti{w}$ in $\Waf$.
\end{rem}

%

\end{document}